\documentclass[11pt]{amsart}

\usepackage[T1]{fontenc}
\usepackage{lmodern}
\usepackage[margin=1in]{geometry}
\usepackage{amsmath,amssymb,amsthm,mathtools}
\usepackage{booktabs,array,tabularx,enumitem,longtable}
\usepackage{microtype}
\usepackage{xurl}
\usepackage{tikz-cd}
\usetikzlibrary{shapes.geometric}
\usepackage{comment}
\usepackage[hidelinks]{hyperref}
\usepackage[all,cmtip]{xy}

\newtheorem{theorem}{Theorem}[section]
\newtheorem{proposition}[theorem]{Proposition}
\newtheorem{lemma}[theorem]{Lemma}

\theoremstyle{definition}

\theoremstyle{remark}

\numberwithin{equation}{section}

\newcommand{\Q}{\mathbb Q}
\newcommand{\C}{\mathbb C}
\newcommand{\Z}{\mathbb Z}
\newcommand{\F}{\mathbb F}
\newcommand{\PP}{\mathbb P}
\newcommand{\gon}{\operatorname{gon}}
\newcommand{\red}{\operatorname{red}}

\newcommand{\Frob}{\operatorname{Frob}}

\newcommand{\End}{\operatorname{End}}

\newcolumntype{Y}{>{\raggedright\arraybackslash}X}

\title{Modular Curves $X_0(N)$ with gonality $7$}

\begin{document}

\begin{abstract}
In this paper we determine all modular curves $X_0(N)$ with $\Q$-gonality equal to $7$. It turns out that there is only one such curve, namely $X_0(133)$. We also fix two previously incorrect results about the gonality of $X_0(N)$.
\end{abstract}
\author{\sc Petar Orli\'c}
\address{Petar Orli\'c \\
University of Zagreb\\  
Bijeni\v{c}ka Cesta 30 \\
10000 Zagreb\\
Croatia}
\email{petar.orlic@math.hr}


\maketitle

\section{Introduction}\label{sec:introduction}

Let $C$ be a smooth projective curve over a field $k$. The $k$-gonality of $C$, denoted by $\textup{gon}_k C$, is the least degree of a non-constant $k$-rational morphism $f:C\to\mathbb{P}^1$.

Gonality is an important invariant that arises in various contexts. It can be understood as a measure of the geometric complexity of a curve and is closely connected to the study of small degree points on it. This is particularly true for modular curves because they are moduli spaces of elliptic curves. For example, the existence of a degree $d$ rational morphism $C\to\PP^1$, together with Hilbert's Irreducibility Theorem \cite[Proposition 3.3.5(2) and Theorem 3.4.1]{serre2016topics}, implies the existence of infinitely many degree $d$ points over $\Q$ on $C$. 

In this paper we focus on the classical modular curve $X_0(N)$ whose properties have been extensively studied. Ogg \cite{Ogg74} determined all hyperelliptic, Bars \cite{Bars99} determined all bielliptic, Hasegawa and Shimura \cite{HasegawaShimura_trig} determined all trigonal curves $X_0(N)$ over $\C$ and $\Q$, and Jeon and Park \cite{JeonPark05} determined all tetragonal curves $X_0(N)$ over $\C$. The paper \cite{NajmanOrlic2023}, which was one of the motivations for this paper, was a follow-up to these results and determined all curves $X_0(N)$ with $\Q$-gonality equal to $4,5,$ or $6$. It also determined the $\Q$-gonality of all curves $X_0(N)$ for $N\leq144$ and for many of those curves also determined their $\C$-gonality.

The AI recently found a few errors in the previously mentioned paper \cite{NajmanOrlic2023}. More specifically, it turns out that the $\Q$-gonality proofs for levels $N=133,145$ contained an error. More details are provided in Section \ref{sec:fixing_mistakes}. This paper fixes these errors and establishes correct $\Q$-gonalities of these two curves.

Another motivation for this paper is the LMFDB \cite{LMFDB}, a database of elliptic curves and modular forms. This database, which is still under development, also contains various data on modular curves, and these results would contribute to it.

The main result of this paper is the following theorem, which is a continuation of $\Q$-gonality results from \cite{NajmanOrlic2023}.

\begin{theorem}\label{thm_main}
    The modular curve $X_0(N)$ has $\Q$-gonality equal to $7$ if and only if $N=133$.
\end{theorem}

Interestingly, this is a very similar situation to the case $\gon_\Q(X_0(N))=5$. There is only one curve $X_0(N)$ with $\Q$-gonality equal to $5$, namely $N=109$. This is an opposite situation to the even cases where there are many curves $X_0(N)$ with that gonality (which is to be expected because Atkin-Lehner quotient maps are of even degree, as well as most degeneracy maps $X_0(N)\to X_0(N)$ because they have degree $\frac{\psi(N)}{\psi(M)}$).

This paper is organized as follows. \begin{itemize}
    \item In Section \ref{sec:gonality_bounds} we provide lower bounds on the $\Q$-gonality of curves $X_0(N)$, either via $\F_p$ gonality or the Castelnuovo-Severi inequality.
    \item In Section \ref{sec:fixing_mistakes} we give proofs for levels $N=133,145,163$. These are the only levels left after Section \ref{sec:gonality_bounds} and the proofs for them are more involved. The proofs for these levels complete the proof of Theorem \ref{thm_main}.
\end{itemize}

\subsection{Data availability and reproducibility statement}
The results in this paper depend on computations carried out in \texttt{Magma} V2.30-1 \cite{magma} and Sagemath 10.7 \cite{sagemath}.

The code for verifying all computations in this paper can be found on the GitHub repository \begin{center}
\url{https://github.com/orlic1/X0-N-_with_Q-gonality_7/tree/main}
\end{center}

The computations were performed on a server at the University of Zagreb equipped with an AMD EPYC 9175F processor (16 cores, up to 4.2 GHz) and 384 GB of RAM.

\subsection{AI disclosure}
Anthropic’s Claude Opus 5.5 reviewed the paper \cite{NajmanOrlic2023}, reran the computations, checked the proofs, and suggested the fixes for levels $N=133,145$ and the argument for level $N=163$. The author independently checked all AI-assisted material and takes full responsibility for the content of the paper.

\section*{Acknowledgments}

Many thanks to Filip Najman for providing the AI-generated materials which were a first step in the creation of this paper.

Maarten Derickx's \texttt{Magma} script was used in some gonality computations. It can be found on \begin{center}
    \url{https://github.com/nt-lib/gonality_bounds/tree/main}
\end{center}

\subsection*{Funding}
The author was supported by the Croatian Science Foundation under the project no. IP-2022-10-5008, by the project ``Implementation of cutting-edge research and its application as part of the Scientific Center of Excellence for Quantum and Complex Systems, and Representations of Lie Algebras", PK.1.1.10.0004, European Union, European Regional Development Fund - Competitiveness and
Cohesion Programme 2021–2027, and by the European Union – NextGenerationEU through the National
Recovery and Resilience Plan 2021–2026, via an institutional grant of the University of Zagreb Faculty of
Science, IK IA 1.1.3, Impact4Math.

\section{Gonality bounds}\label{sec:gonality_bounds}

In this section we will provide results needed for the determination of curves $X_0(N)$ with $\Q$-gonality equal to $7$. First we will give an upper bound on $N$.

\begin{lemma}[Ogg]\cite[Theorem 3.1]{Ogg74}\label{lemmaogg}
    Let $p$ be a prime not dividing $N$. Then the number of $\F_{p^2}$ points on $X_0(N)$ is at least
    $$L_p(N):=\frac{p-1}{12}\psi(N)+2^{\omega(N)}.$$
    Here $\psi(N)$ is the index of the congruence subgroup $\Gamma_0(N)$, equal to $N\prod_{q\mid N}(1+\frac{1}{q})$, and $\omega(N)$ is the number of different prime divisors of $N$.
\end{lemma}
\begin{lemma}\cite[Lemma 3.5]{NajmanOrlic2023}\label{Fp2points}
    Let $C$ be a curve, $p$ a prime of good reduction for $C$, and $q$ a power of $p$. Suppose $\#C(\F_q)>d(q+1)$. Then $\textup{gon}_\Q(C)>d$.
\end{lemma}

Using Ogg's inequality we can easily obtain (see \cite[Lemma 3.2]{HasegawaShimura_trig}, for example) that if $\gon_\Q(X_0(N))\leq7$, then we must have $N\leq396$. Moreover, \cite[Tables 1,2,3]{NajmanOrlic2023} contain results on gonality of curves $X_0(N)$ for $N\leq169$ and many other higher levels $N$. None of the levels $N$ in those tables have $\Q$-gonality equal to $7$. 

After applying Ogg's inequality we are left with $89$ levels $201\leq N\leq 396$. We can eliminate more levels by directly counting $\F_{q}$ points.

\begin{proposition}\label{prop_point_counting}
    The $\Q$-gonality of $X_0(N)$ is at least $8$ for the following levels $N$:
    
    \begin{center}
\begin{longtable}{|c|c|c||c|c|c||c|c|c|}
\caption{Number of $\F_{p^2}$ points on curves $X_0(N)$.}
\label{fpn_points_table}\\
\hline
$N$ & $p$ & $\#X_0(N)(\F_{p^2})$ & $N$ & $p$ & $\#X_0(N)(\F_{p^2})$ & $N$ & $p$ & $\#X_0(N)(\F_{p^2})$

\\
		\hline \hline
$214$ & $3$ & $72$ & $224$ & $3$ & $80$ & $225$ & $2$ & $38$\\
$226$ & $3$ & $72$ & $234$ & $5$ & $196$ & $236$ & $3$ & $78$\\
$237$ & $2$ & $36$ & $242$ & $3$ & $80$ & $246$ & $5$ & $192$\\
$248$ & $3$ & $72$ & $249$ & $2$ & $36$ & $254$ & $3$ & $72$\\
$256$ & $3$ & $80$ & $261$ & $2$ & $38$ & $262$ & $3$ & $80$\\
$265$ & $2$ & $36$ & $270$ & $7$ & $384$ & $301$ & $2$ & $38$\\
$305$ & $2$ & $40$ & $319$ & $2$ & $38$ & $323$ & $2$ & $42$\\
$331$ & $2$ & $37$ & $343$ & $2$ & $36$ & $347$ & $2$ & $37$\\
$349$ & $2$ & $37$ & $353$ & $2$ & $36$ & $361$ & $2$ & $39$\\
$367$ & $3$ & $74$ & $373$ & $2$ & $37$ & $379$ & $2$ & $41$\\
$383$ & $2$ & $38$ & $389$ & $2$ & $37$ & & &\\
		\hline
\end{longtable}
\end{center}
\end{proposition}

\begin{proof}
    The Jacobian $J_0(N)$ has a $\mathbb{Q}$-decomposition into simple abelian varieties
\begin{equation*}\label{eqJac1}
J_0(N)\sim_{\Q} A_{f_1}^{n_1}\times\ldots\times A_{f_k}^{n_k},
\end{equation*}
where $f_i$'s are normalized weight newforms of trivial character and level $d|N$. Let $K_i=\mathbb{Q}(\{a_n(f_i)\})$ be the field of definition of $f_i$. Then we have $\dim A_{f_i}=d_i:=[K_i:\Q]$. The newforms appearing in the decomposition are widely known and available on LMFDB \cite{LMFDB}.

The characteristic polynomial of $\Frob_p$ acting on the Tate module of $A_f$ is
\[\prod_{\sigma:K_i\hookrightarrow\overline{\Q}}(x^2-a_p(f^{\sigma})x+p),\]
where $\sigma$ runs over the set of all $\Q-$embeddings of $K_i$ into a fixed algebraic closure of $\Q$, and $a_p(f^{\sigma})$ denotes the $p$-th coefficient of the $q$-expansion of the modular form $\sigma(f)$.
Thus, the characteristic polynomial of $\Frob_p$ acting on the Tate module of $J_0(N)$) is
\[R_p(x)=\prod_{i=1}^k \left(\prod_{\sigma:K_i\hookrightarrow\overline{\Q}}(x^2-a_p(f^{\sigma})x+p)\right)^{n_i}=\prod_{j=1}^{2g}(x-\alpha_j).\]
Moreover, by the Eichler-Shimura congruence relation in $\End(J_0(N)_{\F_p})$ we have \[T_p=\Frob_p+p\Frob_p^{-1}.\] Therefore, if $T_{p,2}(x)$ is the characteristic polynomial of the Hecke operator $T_p$, then \[R_p(x)=x^{\deg T_{p,2}}\cdot T_{p,2}\left(x+\frac{p}{x}\right).\]

Now, when $p\nmid N$, by the trace formula of Deligne \cite[Exposé III]{deligne77} we have that
\[\#(X_0(N) (\F_{p^n}))=p^n+1-\sum_{j=1}^{2g} \alpha_j^n.\] Notice that $\sum_{j=1}^{2g} \alpha_j^n$ is equal to the trace of the $n$-th power of the companion matrix of $R_p(x)$. This result is important because it gives us a way to compute the number of points over $\F_{p^n}$.

The results in the table were computed using Sage, more details can be found on GitHub. It is now easy to check that in all these case we have $\#X_0(N)(\F_{p^2})>7(p^2+1)$, therefore we can use Lemma \ref{Fp2points} to prove that $\gon_\Q(X_0(N))>7$.
\end{proof}

After this proposition and elimination of all levels solved in \cite[Tables 1,2,3]{NajmanOrlic2023} and \cite[Corollary 3.24]{Orlic2023}, we are left with $46$ levels $N$. These are

\[133,145,163,170,172,174,178,180,182,186,187,189,198,\] \[201, 202, 204, 208, 210, 212, 216, 217, 218, 219, 222, 228, 232, 233, 235, 240,\] \[243, 244, 245, 247, 253, 259, 267, 275, 277, 283, 289, 293, 295, 307, 313, 317, 337.\]

For $N=133$ we will find a degree $7$ rational morphism to $\PP^1$ and for the other $45$ levels we will prove that $\gon_\Q(X_0(N))\geq8$. This will give a proof of Theorem \ref{thm_main}.

\subsection{Computing the $\F_p$-gonality}

We use the fact that for a curve $C/\Q$ and a prime $p$ of good reduction we have \[\gon_{\F_p}(C)\leq\gon_\Q(C).\] Moreoer, $\F_p$-gonality is reasonably computable, as we will see in the next results.

\begin{proposition}\label{Fp_gonality}
    We have a lower bound on the $\Q$-gonality of $X_0(N)$ for the following levels $N$:

        \begin{center}
\begin{longtable}{|c|c|c||c|c|c||c|c|c|}
\caption{Lower bounds on $\gon_{\F_p}(X_0(N))$.}
\label{fp_gonality_table}\\
\hline
$N$ & $p$ & $\gon_{\F_p}(X_0(N))\geq$ & $N$ & $p$ & $\gon_{\F_p}(X_0(N))\geq$ & $N$ & $p$ & $\gon_{\F_p}(X_0(N))\geq$\\
		\hline \hline
$133$ & $3$ & $7$ & $145$ & $2$ & $7$ & $170$ & $3$ & $8$\\
$172$ & $3$ & $8$ & $178$ & $3$ & $8$ & $182$ & $5$ & $8$\\
$187$ & $2$ & $8$ & $189$ & $2$ & $8$ & $201$ & $2$ & $8$\\
$202$ & $3$ & $8$ & $208$ & $3$ & $8$ & $212$ & $3$ & $8$\\
$217$ & $2$ & $8$ & $219$ & $2$ & $8$ & $232$ & $3$ & $8$\\
$233$ & $2$ & $8$ & $235$ & $2$ & $8$ & $243$ & $2$ & $8$\\
$245$ & $2$ & $8$ & $247$ & $2$ & $8$ & $253$ & $2$ & $8$\\
$277$ & $2$ & $8$ & $283$ & $2$ & $8$ & $289$ & $2$ & $8$\\
$307$ & $2$ & $8$ & $313$ & $2$ & $8$ & $317$ & $2$ & $8$\\
$337$ & $2$ & $8$ & & & & & &\\
		\hline
\end{longtable}
\end{center}
\end{proposition}

\begin{proof}
    Using \texttt{Magma}, we compute that there are no functions of degree $\leq d-1$ in $\F_p(X_0(N))$. Let us explain the idea used in the code.

    Let $D\geq0$ be a rational divisor on $X_0(N)$ and $g$ be the genus. By Riemann-Roch theorem we know that \[h^0(D)\geq2\iff h^0(K-D)\geq g-\deg D+1.\] Since $h^0(K)=g$, the problem reduces to computing the dimension of the space of holomorphic differentials vanishing along $D$.

    Choose a basis of holomorphic differentials $\omega_1,\ldots,\omega_g$. At each place $P$ the algorithm computes their local power series expansions. If $u$ is a uniformizing parameter, these expansions have the form \[\omega_j=(a_{j,0}+a_{j,1}u+a_{j,2}u^2+\ldots)du.\] For a divisor $D$, the vanishing conditions on these differentials are encoded in a matrix $M_D$ constructed by combining the relevant expansion coefficients. Its kernel is the space of holomorphic differentials vanishing along $D$. Therefore, $h^0(K-D)=g-\textup{rank}(M_D)$, meaning that \[h^0(D)\geq2\iff \textup{rank}(M_D)\leq \deg(D)-1.\] 
    Moreover, by \cite[Lemma 3.1]{NajmanOrlic2023} we can restrict ourselves to the divisors containing at least \[\Bigg\lceil\frac{\#X_0(N)(\F_p)}{p+1}\Bigg\rceil\] rational points. 
\end{proof}

The computation time for most levels was a few minutes. However, some levels with higher genus took considerably more time. For example, $g(X_0(182))=25$ and $g(X_0(337)=27$. The computation times for these curves were around $19$ days and $30$ days, respectively.

For the next few levels we were able to get a satisfactory lower bound by considering the quotient curve. This drastically reduced the computation time in these cases. Moreover, the genera of these curves $X_0(N)$ are relatively high ($g(X_0(228))=35$, for instance) and a direct computation with $X_0(N)$ would not be feasible in terms of memory and time consumption.

\begin{proposition}\label{Fp_gonality_quotient}
    We have a lower bound on the $\Q$-gonality of $X_0(N)$ for the following levels $N$:

        \begin{center}
\begin{longtable}{|c|c|c|c||c|c|c|c|}
\caption{Lower bounds on $\gon_{\F_p}(X_0(N)/w_d)$.}
\label{fp_gonality_quo_table}\\
\hline
$N$ & $d$ & $p$ & $\gon_{\F_p}(X_0(N)/w_d)\geq$ & $N$ & $d$ & $p$ & $\gon_{\F_p}(X_0(N)/w_d)\geq$\\
		\hline \hline
$180$ & $5$ & $7$ & $8$ & $198$ & $2$ & $5$ & $8$\\
$216$ & $27$ & $5$ & $8$ & $218$ & $2$ & $3$ & $8$\\
$228$ & $76$ & $5$ & $8$ & $244$ & $61$ & $5$ & $8$\\

		\hline
\end{longtable}
\end{center}
\end{proposition}

\begin{proof}
    Similarly to Proposition \ref{Fp_gonality}, we use \texttt{Magma} to compute the lower bound on $\gon_{\F_p}(X_0(N)/w_d)$. After that, we can use \cite[Proposition A.1 (vii)]{Poonen:Gonality}. It says that if $X,Y$ are curves over a field $k$ and $\pi:X\to Y$ a non-constant morphism defined over $k$, then \[\gon_k X\geq \gon_k Y.\] We apply this result with the degree $2$ quotient map $X_0(N)\to X_0(N)/w_d$.
\end{proof}

\subsection{Castelnuovo-Severi inequality}

\begin{proposition}[Castelnuovo-Severi inequality]\label{tm:CS}
Let $k$ be a perfect field, and let $X,\ Y, \ Z$ be curves over $k$. Let non-constant morphisms $\pi_Y:X\rightarrow Y$ and $\pi_Z:X\rightarrow Z$ over $k$ be given, and let their degrees be $m$ and $n$, respectively. Assume that there is no morphism $X\rightarrow X'$ of degree $>1$ through which both $\pi_Y$ and $\pi_Z$ factor. Then the following inequality holds:
\begin{equation} \label{eq:CS}
g(X)\leq m \cdot g(Y)+n\cdot g(Z) +(m-1)(n-1).
\end{equation}
\end{proposition}

Since $\C$ and $\Q$ are both perfect fields, we can use the Castelnuovo-Severi inequality to get lower bounds on both $\C$ and $\Q$-gonalities.

\begin{proposition}\label{CSprop}
The $\C$-gonality of the curve $X_0(N)$ is at least $8$ for the following values of $N$. Here $g$ denotes the genus of the curve $X_0(N)$ and $g'$ denotes the genus of the quotient curve $X_0(N)/w_d$.
\begin{center}
\begin{longtable}{|c|c|c||c|c|c||c|c|c|}

\hline
\addtocounter{table}{-1}
$(N,d)$ & $g$ & $g'$ & $(N,d)$ & $g$ & $g'$ & $(N,d)$ & $g$ & $g'$\\
    \hline
    
$(174,87)$ & $27$ & $8$ & $(186,62)$ & $29$ & $11$ & $(204,68)$ & $31$ & $12$\\
$(210,14)$ & $41$ & $17$ & $(222,111)$ & $35$ & $10$ & $(240,80)$ & $37$ & $15$\\
$(259,259)$ & $23$ & $8$ & $(267,89)$ & $29$ & $9$ & $(275,275)$ & $25$ & $9$\\
$(293,293)$ & $24$ & $8$ & $(295,295)$ & $29$ & $11$ & & &\\
 
    \hline

\caption*{}
\end{longtable}
\end{center}
\end{proposition}

\begin{proof}
    Suppose that there exists a morphism $f:X_0(N)/w_d\to\PP^1$ of degree $\deg f\leq7$. We apply the Castelnuovo-Severi inequality with $f$ and a degree $2$ quotient map $\pi: X_0(N)/w_d\to X_0(N)/\left<w_d,w_{d'}\right>$. In all cases here we have \[g(X_0(N))\geq2g(X_0(N)/w_d)+7.\] Therefore $f$ must factor through $\pi$ (because $\deg\pi=2$ is a prime). This is impossible if $2\nmid \deg f$, and if $2\mid \deg f$ we would have $\gon_\C(X_0(N)/w_d)\leq\frac{\deg f}{2}\leq3$. However, the results of \cite{FurumotoHasegawa1999} and \cite{HasegawaShimura1999} tell us that these curves $X_0(N)/w_d$ are neither hyperelliptic nor trigonal over $\C$ and we get a contradiction.
\end{proof}

\section{The arguments for $N=133,145,163$}\label{sec:fixing_mistakes}

\subsection{$N=133$}

We will prove that the $\Q$-gonality of $X_0(133)$ is equal to $7$, not $8$ as stated in \cite{NajmanOrlic2023}. The Mordell-Weil sieving \texttt{Magma} code used there has an error; it erroneously used the Atkin-Lehner involution $w_7$ instead of the Fricke involution $w_{133}$, and hence obtained the wrong result.

First, Proposition \ref{Fp_gonality} tells us that $\gon_\Q(X_0(133))\geq7$. It remains to find a degree $7$ rational morphism to $\PP^1$.

We know that $X_0(133)$ has genus $11$ and $4$ rational cusps, these are the only rational points on that curve. We construct a canonical model of this curve inside $\PP^{10}$ where all Atkin-Lehner involutions act as $\pm1$ on coordinates (the same model as in \cite{NajmanOrlic2023}). In that model the cusps are \begin{align*}
    C_1&=[1,0,1,0,0,0,1,0,1,0,0],\\
    C_2&=[-1,0,-1,0,0,0,1,0,1,0,0],\\
    C_3&=[-1,0,1,0,0,0,-1,0,1,0,0],\\
    C_4&=[1,0,-1,0,0,0,-1,0,1,0,0].
\end{align*}

The quotient curve $X_0^+(133)$ has genus $5$. The forms $x_6,x_5,x_8-x_4$ map $X_0(133)$ into a plane septic model of $X_0^+(133)$. This model has a node at $(1:1:0)$. Pulling this node back to $X_0(133)$ and removing
the cusp base locus gives a degree $4$ rational divisor $E$, and the code checked that it is irreducible.

The code constructs $E$ from the linear forms vanishing on $E$ (that is, $H^0(K-E)$). Riemann-Roch theorem gives us that \[h^0(K-E)=h^0(E)-(\deg E+1-g)=h^0(E)+6,\] and we have $h^0(E)=1$ (otherwise we would have a contradiction with $\gon_\Q\geq7$). Thus, there are $7$ such forms.

We now combine $E$ with the rational points into a degree $7$ effective divisor. It turns out that the divisor \[D=E+2C_1+C_3\] has Riemann-Roch dimension equal to $2$ (the code computes that $h^0(K-D)=5$, and the statement follows from Riemann-Roch), thus we get $\gon_\Q(X_0(133))=7$.

\subsection{$N=145$}

We will prove that the $\Q$-gonality of $X_0(145)$ is equal to $8$. The Mordell-Weil sieving \texttt{Magma} code used in \cite{NajmanOrlic2023} did not go through all possible degree $6$ divisors, yielding an invalid proof. Actually, it turns out that going through all divisors yields a set $\{0\}\subseteq J_0(145)^-$, instead of an empty set which would disprove the existence of degree $6$ maps to $\PP^1$. We will explain more in the proof that $\gon_\Q=8$.

The curve $X_0(145)$ has genus $13$ and $4$ rational cusps, these are the only rational points on that curve \cite{mazur1978rational}. The quotient curve $X_0^+(145)$ has genus $5$. Proposition \ref{Fp_gonality} tells us that $\gon_\Q(X_0(145))\geq7$. In addition, the curve $X_0(145)/w_{29}$ is of genus $4$ (and has a rational cusp), so we can construct a degree $8$ rational map as a composition \[X_0(145)\to X_0(145)/w_{29}\to \PP^1.\] It remains to show that $\gon_\Q\neq7$.

Let $X:=X_0(145), Y:=X_0^+(145), w:=w_{145}$. Riemann-Hurwitz formula applied to the degree $2$ quotient map $\pi:X\to Y$ tells us that there are $8$ ramification points. It turns out that the ramification divisor $R$ is a single reduced closed point of degree $8$ and $h^0(R)=2$.

Take a nonzero rational function $u$ on $X$ with $u \circ w = -u$. The ratio of a minus cusp form to a plus cusp form works and is defined over $\Q$. At each fixed point $u$ has odd order, and away from them its zeros and poles come in $w$-symmetric pairs. So \[\textup{div } u=\pi^*U+R\] for some rational divisor $U$ on $Y$ of degree $-4$. 
By splitting functions on $X$ into even and odd parts we get \[\pi_*\mathcal{O}_X=\mathcal{O}_Y\oplus \mathcal{O}_Y(U).\]
We have $\mathcal{O}_X(R)=\pi^*(-U)$ and \[2=h^0(\mathcal{O}_X(R))=h^0(\pi^*(-U))=h^0(-U)+h^0(\mathcal{O}_Y)=h^0(-U)+1\implies h^0(-U)=1.\]

One can easily check on LMFDB \cite{LMFDB} by looking at newforms of level dividing $145$ that $J_0(145)(\Q)^-$ is of rank $0$ and using the result of Kolyvagin and Logachev \cite{KolyvaginLogachev89}. We define the map \[\mu:\textup{Pic}^d(X)\to J_0(N), \qquad \mu(D):= D-w(D).\] We obviously have that $\mu(J_0(145))\subseteq J_0(145)^-$.

Suppose that $M\geq0$ is a degree $7$ rational divisor with $h^0(M)\geq2$. Consider a divisor \[\mu(M)=M-w(M)=(1-w)(M-7C)+7(C-w(C))\in J_0(145)_{\textup{tors}}.\] ($\mu(M)=M\otimes (w^*M)^{-1}$ in the language of line bundles). This inclusion holds because $J_0(145)(\Q)^-\subseteq J_0(145)_{\textup{tors}}$ and $C-w(C)$ is in a cuspidal subgroup which is contained in the torsion. Since $h^0(M)\geq2$, we may assume that $M=C+E$ for some fixed rational point $C$ and $E\geq0$ of degree $6$ (instead of $f\in H^0(M)$ we can take $\frac{1}{f-f(C)}$).

The Mordell-Weil sieve we use now is similar to the one in \cite[Proposition 6.1]{NajmanOrlic2023}, but we use $\mu$ instead of $2\mu$. Suppose $p>2$ is a prime of good reduction for $X$. Then we have the following commutative diagram
$$
\xymatrix{
& W_d^1(X)(\Q) \ar[r]^{\mu} \ar[d]^{\red_p} & J_0(N)(\Q)_{\textup{tors}} \ar[d]^{\red_p} \\
& W_d^1(X)(\F_p) \ar[r]^{\mu} & J_0(N)(\F_p)}, 
$$
where the vertical maps are reduction modulo $p$. If there is some $D\in W_d^1(X)(\Q)$, then \[\mu(D)\in\red_p^{-1}(\mu(W_d^1(X)(\F_p))).\] The set $W_d^1(X)(\F_p)$ can be computed in practice by simply finding all the effective degree $d$ divisors whose Riemann-Roch spaces have dimension $\geq 2$. For $N=145$ and $p=3$ there are $5750$ effective degree $6$ divisors $E$ and exactly $7$ divisors $C+E$ have $h^0\geq2$ over $\F_3$. From the code we get that their $\mu$-images are exactly the $\mu$-images of the $4$ rational cusps.

Note that $\red_p$ is injective on the torsion of $J_0(N)(\Q)$ \cite[Appendix]{katz81}. Hence we have $\mu(M)\sim\mu(P)$ over $\Q$ for some rational cusp $P$.

Define $A:=M-P$. Then $A$ is defined over $\Q$ and $\mu(A)=0$. We also have that \[h^0(A)=h^0(M-P)\geq h^0(M)-1\geq2-1=1.\] Moreover, we must have $h^0(A)=1$ because $\deg A=6$ and $\gon_\Q X\geq7$. Thus there is a unique rational divisor $E\geq0$ equivalent to $A$. Uniqueness implies $w(E)=E$ (because $w(A)\sim A$ and $\gon_\Q X\geq7$). Since $R$ is a Galois orbit of $8$ points, $E$ must be unramified and we have $\mathcal{O}_X(E)=\pi^*D$ for a rational divisor $D$ with $\deg D=3$ and $h^0(D)=1$ (the argument for $h^0=1$ is the same as with $R$ and $-U$).

Let $Q=\pi(P)\in Y$. Since $P$ is unramified, we have $\pi^*Q=P+w(P)$. Define a map \[\Phi:H^0(E+P)\to \Q(X), \qquad \Phi(f)=f-f\circ w.\]

\begin{lemma}\label{lem_ker_phi_145}
    $\ker\Phi=H^0(E)\cong \Q$.
\end{lemma}
\begin{proof}
    If $f=f\circ w$ for some $f\in H^0(E+P)$, then $\textup{div } f+E+P\geq0$. Applying $w$ gives \[\textup{div } f=\textup{div} (f\circ w)\geq -E-w(P).\] Therefore, $\textup{div } f\geq -E$ because $P\neq w(P)$. On the other hand, $h^0(E)=h^0(A)=1$ and $E\geq0$, meaning that $H^0(E)$ consists only of constant functions, and they are obviously invariant under $w$.
\end{proof}

\begin{lemma}\label{lem_im_phi_145}
    $\dim \textup{im } \Phi\leq h^0(D+Q+U)$.
\end{lemma}
\begin{proof}
    The image of $\Phi$ lies in the anti-invariant part of $H^0(\pi^*(D+Q))$ because \[\textup{div } f\geq -E-P\implies \textup{div} (f\circ w)\geq-E-w(P)\implies \textup{div} (\Phi(f))\geq -E-P-w(P)=-\pi^*(D+Q).\]
    We can write $\Phi(f)=u\cdot(h\circ \pi)$, which is possible because $\frac{\Phi(f)}{u}$ is invariant. We now have \[\textup{div} (h\circ\pi)\geq -\pi^*(D+Q+U) - R.\] Let us check this condition at every point $y\in Y$. Consider an integer $x:=\textup{ord}_y h+\textup{mult}_y (D+Q+U)$. We have $2$ cases. \begin{itemize}
        \item $y$ is not a branch point: Then \[\textup{div}(\Phi(f))=\textup{div }u+\textup{div}(h\circ\pi)=\pi^*U+R+\pi^*\textup{div } h \implies \textup{mult}_y U+\textup{ord}_y h\geq-\textup{mult}_y (D+Q).\] Hence $x\geq0$.
        \item $y$ is a branch point: Then \[\textup{div}(\Phi(f))=\textup{div }u+\textup{div}(h\circ\pi)=\pi^*U+R+\pi^*\textup{div } h \implies 2\textup{mult}_y U+1+2\textup{ord}_y h\geq-2\textup{mult}_y (D+Q).\] Therefore, $2x+1\geq0$, so $x\geq0$ because it is an integer.
    \end{itemize}
    In other words, we have proven that $\textup{div } h+D+Q+U\geq0$, or equivalently $h\in H^0(D+Q+U)$. Since composition by $\pi$ is injective in $\Q(X)$, this implies that $\dim \textup{im } \Phi\leq h^0(D+Q+U)$, as desired.
\end{proof}

Now we have \[2\leq h^0(M)=h^0(E+P)=\dim\ker\Phi+\dim\textup{im }\Phi\leq1+h^0(D+Q+U)\implies h^0(D+Q+U)\geq1.\] However, $\deg (D+Q+U)=3+1-4=0$, meaning that $D+Q+U$ is principal and \[\pi^*(D+Q+U)=E+P+w(P)+\pi^*(U)\sim 0\implies A+P+w(P)-R\sim 0\implies M\sim R-w(P).\] On the other hand, since the constant function $1$ is in $H^0(R)\setminus H^0(R-w(P))$ (because $P$ is a cusp and we computed that $R$ is a closed degree $8$ point), this implies that $h^0(M)\leq h^0(R)-1=1$, a contradiction with our assumption. This proves that $\gon_\Q(X)\geq8$.

\subsection{$N=163$}

We will prove that the $\Q$-gonality of $X_0(163)$ is equal to $8$. This curve has genus $13$ and $3$ rational points: cusps $0,\infty$ and $1$ noncuspidal point \cite{mazur1978rational}. The genus $6$ quotient $X_0^+(163)$ is tetragonal over $\Q$ by \cite{Orlic2024} and the composition map \[X_0(163)\to X_0^+(163)\to\PP^1\] provides a rational degree $8$ morphism to $\PP^1$. It remains to prove that $\gon_\Q\geq8$.

Let $X:=X_0(163), Y:=X_0^+(163), w:=w_{163}$. Riemann-Hurwitz formula applied to the quotient map $\pi:X\to Y$ tells us that there are $4$ fixed points of $w$ on $X$. These are CM points. It turns out that we have a rational point $R_0$ and a closed point $R_3$ of degree $3$. Define $R:=R_0+R_3$ and write $Q_0:=\pi(R_0), B_3:=\pi(R_3)$.

Suppose that $M\geq0$ is a degree $7$ rational divisor with $h^0(M)\geq2$ (if there is a function of degree $d\leq6$, we can just add $7-d$ cusps to its polar divisor to get $M$). We may assume that $M=C+E$ for some rational point $C$ and $E\geq0$ of degree $6$.

We define $U$ in the same way as for $N=145$. Then $U$ is rational, $\deg U=-2$, and $\mathcal{O}_X(R)=\pi^*(-U)$. 

Now we will do the Mordell-Weil sieve, similarly as for $N=145$ and with the same $\mu$. One can easily check on LMFDB by looking at newforms of level $163$ that $J_0(163)(\Q)^-$ is of rank $0$. Since $163$ is prime, the torsion of $J_0(163)$ is the cuspidal subgroup \[J_0(163)_{\textup{tors}}=\left<[0-\infty]\right>\cong \Z/27\Z.\] We take $p=3$. The code finds that there are $17514$ divisors $E$ over $\F_3$ and exactly $3$ of the divisors $C+E$ have $h^0\geq2$. Of their $\mu$-images, only one lies in the order-$27$ cuspidal subgroup, and it is $0\in\Z/27\Z$. Since the reduction mod $p$ is injective on the torsion, this means that $\mu(M)=0\implies w(M)\sim M$. Thus there is a rational function $g$ with $\textup{div }g=w(M)-M$.

Now we have \[\textup{div } g\cdot(g\circ w)=w(M)-M+M-w(M)=0 \implies g\cdot(g\circ w)=c\in\Q^\times.\] However, when we plug in $R_0$, we get $w(R_0)=R_0$ and $g(R_0)^2=c\neq0$. If we replace $g$ by $\frac{g}{g(R_0)}$, we get $g\cdot(g\circ w)=1$. Let $k=1+g$. Then $k\neq0$ (since $g(R_0)=1$) and \[\frac{k}{k\circ w}=\frac{1+g}{1+\frac{1}{g}}=g.\] Set $E=M+\textup{div }k$. Then \[w(E)=w(M)+\textup{div }(k\circ w)=M+\textup{div }g+\textup{div }(k\circ w)=M+\textup{div }k=E.\] So $E$ is $w$-invariant of degree $7$, defined over $\Q$, and $E\sim M\implies h^0(E)=h^0(M)\geq2$ by assumption ($E$ is not necessarily effective, however).

Since $E$ is $w$-invariant, non-fixed points of $E$ come in pairs $P+w(P)=\pi^*(\pi(P))$. Thus we have $E=\pi^*F+S$, where $S$ is a reduced set of fixed points. Since $E$ is rational and $\deg E$ is odd, this means that $S=R_0$ or $S=R_3$. However, we have \[\textup{div }u=\pi^*U+R_0+R_3\] (the definition of $u$ and its properties are the same as for $N=145$, the only difference is that here $\deg U=-2$). Therefore, we may assume that \[E=\pi^*F+R_0, \qquad \deg F=3.\]

\begin{lemma}\label{lem_dim_l_e_163}
    $h^0(E)=h^0(F)+h^0(F+U+Q_0)$.
\end{lemma}
\begin{proof}
    The map $f\mapsto f\circ w$ is invariant on $H^0(E)$ and we have \[f=\frac{f+f\circ w}{2}+\frac{f-f\circ w}{2}.\] Hence, $H^0(E)$ splits into an odd and even part. \begin{itemize}
        \item Even part: $f=h\circ\pi$, at $R_0$ the order of $f$ is even, so the allowed pole of order $1$ does not help. Hence $f\in H^0(\pi^*F)\iff h\in H^0_Y(F)$.
        \item Odd part: $f=u\cdot(h\circ \pi)$, then \[\textup{div }f=\textup{div }u+\textup{div}(h\circ \pi)=\pi^*U+R_0+R_3+\pi^*\textup{div }h\geq -E=-\pi^*F-R_0.\] For every point $y\in Y$ consider an integer $x=\textup{ord}_y h+\textup{mult}_y(F+U)$. We then have $2x+2\cdot1_{R_0}+1_{R_3}\geq0$ at branch points and $x\geq0$ otherwise. Thus $x\geq-1$ at $Q_0$ and $x\geq0$ elsewhere, implying that $h\in H^0_Y(F+U+Q_0)$.
    \end{itemize}
    This proves the desired statement.
\end{proof}

Since $\deg F=3$, $\deg(F+U+Q_0)=3-2+1=2$, and $Y$ is neither hyperelliptic nor trigonal by \cite{FurumotoHasegawa1999, HasegawaShimura1999}, we must have $h^0(F),h^0(F+U+Q_0)\leq1$. By assumption we have $h^0(E)\geq2$, hence $h^0(F)=h^0(F+U+Q_0)=1$. Let $T\geq0$ be the unique divisor with $T\sim F+U+Q_0$. 

\begin{lemma}\label{lem_very_ample_163}
    $A:=K_Y-U-B_3$ is very ample of degree $9$ and $h^0(A)=4$.
\end{lemma}
\begin{proof}
    We have $\deg A=10+2-3=9$. Since \[\pi_*\mathcal{O}_X=\mathcal{O}_Y\oplus \mathcal{O}_Y(U),\] the invariant differentials on $X$ are pullbacks $\pi^*H^0(K_Y)$, and the anti-invariant differentials correspond to $H^0_Y(K_Y-U)$. By Riemann-Roch we know that \[h^0(K_Y-U)=h^0(U)-(\deg U+1-g(Y))=0+2-1+6=7.\] Locally at a ramification point, an anti-invariant differential has the form $f(t^2)dt$, so vanishing on $R_3$ forces vanishing on $2R_3$. Thus, the anti-invariant differentials vanishing on $R_3$ correspond to $H^0(K_Y-U-B_3)=H^0(A)$. The \texttt{Magma} code finds that $h^0(A)=4$ by considering the $7$ newforms of level $163$ with $-$ Fricke sign.
    
    Let $\{\eta_1,\eta_2,\eta_3,\eta_4\}$ be the basis of $H^0(A)$. Let \[\Psi=(\eta_1,\eta_2,\eta_3,\eta_4):X\to\PP^3.\] Each ratio $\frac{\eta_i}{\eta_j}$ is invariant under $w$, hence \[\Psi=\Phi_A\circ \pi, \] where $\Phi_A:Y\to\PP^3$ is the map defined by the complete linear system $|A|$. The common zeros of $\eta_i$ on $X$ are $\pi^*B_3$ plus the pullback of any base points of $|A|$. The \texttt{Magma} code checks that the saturated base ideal is the ideal of $2R_3 = \pi^*B_3$. Therefore $|A|$ has no base points and $\Phi_A$ is a morphism.

    The code checks that $\textup{im }\Psi$ has degree $9$ in $\PP^3$ and arithmetic genus $6=g(Y)$. With this we can conclude that \[18=24-6=\deg K_X-2\deg R_3=\deg\Psi\cdot \deg(\textup{im }\Psi)=2\cdot\deg \Phi_A\cdot9\implies \deg \Phi_A=1,\] so $\Phi_A$ is birational onto its image which is a smooth curve. This proves that $A$ is very ample since a birational morphism between smooth projective curves is an isomorphism.
\end{proof}

Very ampleness, stated with divisors, says that $h^0(A-T)=h^0(A)-\deg T$ for every $T\geq0$ of degree $\deg T\leq2$. If we take our $T\sim F+U+Q_0$, we get \[A-T\sim K_Y-U-B_3-F-U-Q_0=K_Y-F-(2U+Q_0+B_3).\] However, $u^2$ is $w$-invariant, so it factorizes through $\pi$. We have \[u^2=v\circ\pi\implies \textup{div}(u^2)=2\pi^*U+2R_0+2R_3\implies \textup{div }v=2U+Q_0+B_3\implies A-T\sim K_Y-F.\] Riemann-Roch theorem tells us that \[h^0(A-T)=h^0(K_Y-F)=h^0(F)-(\deg F+1-g(Y))=3=h^0(A)-1.\] However, $\deg T=3-2+1=2$ and we get a contradiction. Thus $\gon_\Q X\geq8$.

\bibliographystyle{siam}
\bibliography{biblio}\vspace{0.75in}

\end{document}